\documentclass[11pt]{amsart}

\usepackage{amsmath}
\usepackage{amssymb}
\usepackage{amsthm}
\usepackage{mathtools}
\usepackage[T1]{fontenc}
\usepackage[utf8]{inputenc}

\theoremstyle{plain}
\newtheorem{theorem}{Theorem}[section]
\newtheorem{lemma}[theorem]{Lemma}
\newtheorem{proposition}[theorem]{Proposition}
\newtheorem{corollary}[theorem]{Corollary}
\newtheorem*{theoremA}{Theorem A}
\newtheorem*{theoremB}{Theorem B}
\newtheorem*{theoremC}{Theorem C}
\newtheorem*{theoremD}{Theorem D}

\theoremstyle{definition}
\newtheorem{question}[theorem]{Question}

\newcommand{\Homeo}{\operatorname{Homeo}}
\newcommand{\Isom}{\operatorname{Isom}}
\newcommand{\diam}{\operatorname{diam}}
\newcommand{\dist}{\operatorname{dist}}
\newcommand{\Per}{\operatorname{Per}}
\newcommand{\ext}{\operatorname{ext}}
\newcommand{\Sym}{\operatorname{Sym}}

\title{Minimal sets of countable group actions on regular curves}


\begin{document}
	
\author[E.~Shi]{Enhui Shi}
\address{School of Mathematics and Sciences, Soochow University, Suzhou, Jiangsu 215006, China}
\email{ehshi@suda.edu.cn}

\author[S.~Wang]{Suhua Wang*}
\thanks{*Corresponding author}
\address{School of Mathematical Science, Suzhou Vocational University, Suzhou, Jiangsu
	215104, People’s Republic of China}
\email{wangsuhuajust@163.com}

\author[Y.~Yu]{Yifei Yu}
\address{School of Mathematics and Sciences, Soochow University, Suzhou, Jiangsu 215006, China}
	\email{20254207022@stu.suda.edu.cn}

\begin{abstract}
We study minimal sets of countable groups acting by homeomorphisms on regular
curves. If the ambient action has a finite orbit, every infinite minimal set is an
equicontinuous Cantor system. If there is no finite orbit, the minimal set is unique
and remains minimal for every finite-index subgroup. A minimal set supporting an
invariant probability is a finite orbit, an equicontinuous Cantor system, a circle
isometry action, or a Cantor extension obtained by splitting a countable invariant
subset of such a circle. Without an invariant probability, infinitely many
nondegenerate components can occur. Their homeomorphism types determine the
underlying space, and every nonempty countable family of types is realized. The
two-generator free group admits an example with arc components on a dendrite,
whereas its restricted wreath product with the integers admits none on any regular
curve.
\end{abstract}

\maketitle
\section{Introduction}

\subsection{Background.}
Throughout this paper, $G$ is a countable discrete group and all group actions are
by homeomorphisms. A \emph{continuum} is a nonempty compact connected metric space.
A nondegenerate continuum $X$ is called a \emph{regular curve} if every point has a
neighborhood basis consisting of open sets with finite boundary. This class contains
graphs and dendrites. A \emph{minimal set} of a $G$-action on $X$ is a nonempty
closed invariant subset $M$ such that the orbit $Gx$ is dense in $M$ for every
$x\in M$. We study the topology of these sets and the actions induced on them.

For a homeomorphism of a regular curve, every infinite minimal set is a Cantor set
or a circle; if the homeomorphism has a periodic point, each infinite minimal set is
an adding machine; if it has no periodic point, its minimal set is unique
\cite{ref15,ref3,ref1}. Related results hold for group actions on local dendrites: a
finite ambient orbit implies equicontinuity on every minimal set \cite{ref7}; on a
dendrite, an action without a finite orbit has a unique minimal set \cite{ref6,ref10}.
Minimal sets on a finite graph different from a circle are finite, while those on a
local dendrite containing a circle are finite orbits, Cantor sets or circles
\cite[Theorems 3.6 and 3.9]{ref13}. The relation between invariant measures and free
subgroups has also been studied for actions on curves. Duchesne and Monod proved
alternatives and fixed-point results for actions on dendrites and local dendrites
\cite{ref6}. For a minimal action on a totally regular curve, Shi, Xu and Ye proved
that the action is conjugate to circle isometries or contains a quasi-Schottky
subgroup \cite[Theorem 1.3]{ref18}, which is further extended to non-minimal actions by Xu \cite{ref21}. 

\subsection{Main results.}
We first determine the effect of a finite orbit in the ambient curve.

\begin{theoremA}
Let $G$ act on a regular curve $X$. If $G$ has a finite orbit in $X$, then each
infinite minimal set is an equicontinuous Cantor system. Otherwise $X$ has a unique
minimal set, which remains the unique minimal set for every finite-index subgroup of
$G$. A minimal set with finitely many connected components is either connected or a
finite orbit.
\end{theoremA}

The second result classifies minimal sets supporting an invariant probability. A
\emph{split circle} over a countable invariant subset $D$ is obtained by replacing
each point of $D$ by its two one-sided copies in the circular order. The action lifts
to this space by its action on the order and on the underlying circle. The
construction is given in Subsection~\ref{subsec:split}.

\begin{theoremB}
Let $M$ be a minimal set of a regular-curve $G$-action. If $M$ supports an invariant
Borel probability measure, then $(M,G)$ is conjugate to one of the following systems:
\begin{enumerate}
\item[(1)] a finite transitive $G$-set;
\item[(2)] an equicontinuous minimal $G$-action on a Cantor set;
\item[(3)] a minimal action by circle isometries;
\item[(4)] the lift of a minimal circle isometry action to the circle split over a
nonempty countable invariant set.
\end{enumerate}
Conversely, each system of these four types can be realized as a minimal set of a
regular-curve action of the same group. In (3) and (4), $M$ is the unique minimal set
of the ambient action. In (4), $M$ is a Cantor set, and the collapse map is the
maximal equicontinuous factor map. Its exceptional fibers consist of exactly two
points. In (2), this factor map is a conjugacy and the ambient action has a finite
orbit.
\end{theoremB}

Every action of an amenable group on a compact space admits an invariant probability.
Thus Theorem B applies to every minimal set of an amenable group action on a regular
curve. For actions on dendrites, the infinite case is equicontinuous, as proved by
Shi and Ye \cite{ref19}. For actions on a circle, the split-circle case can occur even
when every element of the acting group has finite order. For a minimal set $M$, let
$q\colon M\to Q$ be the quotient map obtained by collapsing each connected component
to a point. The next theorem describes the case of infinitely many components.

\begin{theoremC}
If $M$ has infinitely many components, $Q$ is a Cantor minimal system, and $q$ is
almost one-to-one with at most countably many nondegenerate fibers. The
homeomorphism type of $M$ is determined by the set of homeomorphism types of its
nondegenerate components. Every nonempty finite or countably infinite set of types of
nondegenerate regular curves is realized in a proper minimal set for the free group of
countably infinite rank. The component stabilizers of the nondegenerate components in
this construction are trivial. If all the prescribed curves are dendrites, the ambient
space may be chosen to be a dendrite.
\end{theoremC}

When there are finitely many nondegenerate component types, the space in Theorem C is
the dense amalgam of \'Swi\k{a}tkowski \cite[Theorem 0.2]{ref20}. Our proof also
applies to countably many types. The classification is a classification of the
underlying compact space; it does not determine the action on that space. The choice
of the acting group imposes further restrictions. Write $F_2$ for the free group on
two generators and $\wr$ for the restricted wreath product.

\begin{theoremD}
If a minimal set $M$ has infinitely many components and at least one is
nondegenerate, the image of $G$ in $\Homeo(M)$ contains a nonabelian free subgroup.
The group $F_2$ itself realizes such a minimal set with arc components on a dendrite.
The group $F_2\wr\mathbb{Z}$ contains $F_2$, but every minimal set of every one of its
regular-curve actions supports an invariant probability measure.
\end{theoremD}

The first assertion also follows from Theorem B and the invariant-measure alternative
in Xu's companion manuscript \cite[Theorem 1.3]{ref21}. We give a proof for minimal
sets which also shows that, in the absence of an invariant probability, the
centralizer of the restricted action is finite. This is used in the wreath-product
argument. The commuting-subgroup obstructions of Duchesne and Monod provide related
results for dendrites and local dendrites \cite[Theorems 1.4 and 11.3]{ref6}.

\subsection{Outline of the proofs.}
The proof of Theorem A starts with open neighborhoods separating an invariant compact
set from a disjoint infinite minimal set. The union of the neighborhoods with least
boundary cardinality has finite invariant boundary. For a connected invariant compact
set, this construction gives an exhaustion of a larger invariant open set by connected
invariant compact sets. The complementary components whose closures meet its frontier
form finite permutation systems. Their inverse limit maps onto the frontier, which
gives the required equicontinuity. For Theorem B, we use Theorems 2.8 and 2.9 to
obtain an almost one-to-one equicontinuous factor. We prove that the resulting almost
one-to-one equicontinuous factor has finite fibers, with only finitely many fibers of
any given positive diameter. Theorem A then distinguishes the Cantor factor from the
circle factor. In the latter case the circular order shows that each fiber has at most
two points. The topological part of Theorem C is proved by successive clopen
partitions around components. For the realization, we replace vertices of a countable
tree by the prescribed regular curves and join them by intervals. An HNN extension
gives the two-generator example in Theorem D. For the free-subgroup assertion, an
irreducible compact affine system gives a finite set toward which every finite subset
can be moved. Finite boundaries then give the inclusions needed for the ping-pong
argument. The same finite-set property yields the centralizer conclusion used for
wreath products.

\subsection{Organization of the paper.}
Section 2 contains the preliminary definitions and results, including the component
quotient. Section 3 proves Theorem A. Section 4 gives the classification with an
invariant probability and its consequences for cyclic subactions. The models in
Section 5 prove the realization assertions in Theorem B. Section 6 treats the topology
and realization of component types. Section 7 proves the free-subgroup alternative and
the wreath-product obstruction. Finally, Section 8 gives consequences for a prescribed
acting group and states the remaining realization question.
\section{Preliminaries}

In this section we recall the notions and results which will be used in the sequel.
The acting groups $G$ are countable and discrete, and all actions are by
homeomorphisms on compact metric spaces.

\subsection{Group actions and invariant measures.}
Let $G$ act on a compact metric space $X$. We write $gx$ for the image of $x$ under
$g\in G$, and $Gx=\{gx:g\in G\}$ for its orbit. A nonempty closed invariant subset
$M$ of $X$ is a minimal set if $\overline{Gx}=M$ for every $x\in M$. A factor map
between two $G$-systems is a continuous equivariant surjection. It is called monotone
if all its fibers are connected. A factor map is almost one-to-one if the points with
singleton fibers are dense in the range. For minimal systems their inverse image is
then dense in the domain. A system $(X,G)$ is equicontinuous if the family of action
maps is equicontinuous. Equivalently, $X$ admits a compatible $G$-invariant metric.
Its maximal equicontinuous factor is the equicontinuous factor through which every
equicontinuous factor factors. A minimal system is almost automorphic if it is an
almost one-to-one extension of a minimal equicontinuous system.

We write $\mathcal{E}(X,G)$ for the closure of the action maps in $X^X$ with the
topology of pointwise convergence; this is the enveloping semigroup. For compact
metric systems we use the equivalent definition that $(X,G)$ is tame if
$|\mathcal{E}(X,G)|\le 2^{\aleph_0}$. A system is null if its topological sequence
entropy vanishes along every sequence in $G$.

Let $\mathcal{P}(X)$ be the space of Borel probability measures on $X$ with the
weak$^*$ topology. It is a compact metrizable convex space. The induced action is
given by $(g_*\mu)(B)=\mu(g^{-1}B)$ for every Borel set $B$. We write
$\mathcal{P}_G(X)$ for the invariant measures and $\delta_x$ for the point mass at
$x$. We also write $g\mu$ for $g_*\mu$. An invariant probability on an infinite
minimal system has full support and no atoms: its support is a nonempty closed
invariant set, and an atom would have a finite orbit.

Two points $x,y\in X$ are proximal if there are $g_n\in G$ with
$d(g_nx,g_ny)\to 0$. The action is proximal if every pair is proximal. It is strongly
proximal if the orbit closure of every member of $\mathcal{P}(X)$ contains a point
mass. It is extremely proximal if every proper closed subset can be moved into every
nonempty open subset. Recall that $G$ is amenable if it admits a left invariant mean
on $\ell^\infty(G)$. We will use the following elementary consequence.

\begin{lemma}
If $G$ is amenable, then $\mathcal{P}_G(X)\neq\emptyset$ for every nonempty compact
metric $G$-space $X$. More generally, every continuous affine action of $G$ on a
nonempty compact metrizable convex subset of a locally convex space has a fixed
point.
\end{lemma}

\begin{proof}
Fix $x\in X$ and let $m$ be a left invariant mean. The functional
$f\mapsto m\bigl(g\mapsto f(gx)\bigr)$ on $C(X)$ is positive, normalized and
$G$-invariant, and hence defines an invariant probability. For an affine action on a
compact convex set, the barycenter of an invariant probability is a fixed point.
\end{proof}

\subsection{Continua and hyperspaces.}
A continuum is a nonempty compact connected metric space. It is nondegenerate if it
has more than one point. A regular curve is a nondegenerate continuum $X$ such that
each point has a neighborhood basis of open sets with finite boundary. A point has
order at most $n$ in $X$ if these neighborhoods can be chosen with at most $n$
boundary points. A dendrite is a locally connected continuum containing no simple
closed curve. A continuum is Suslinian if each family of pairwise disjoint
nondegenerate subcontinua is countable. A compact space is homogeneous if its
homeomorphism group acts transitively.

Regular curves are locally connected, every nondegenerate subcontinuum of a regular
curve is regular, and every nondegenerate dendrite is regular. Connected open subsets
of a locally connected metric continuum are arcwise connected. We use these standard
facts in continuum theory as in \cite{ref14}. The assertion about subcontinua also
follows from $\partial_C(U\cap C)\subset C\cap\partial_X U$ for $U$ open in $X$. The
boundary and closure of a set are taken in the ambient space unless otherwise
indicated.

For a compact metric space $(X,d)$, let $2^X$ be the space of nonempty compact
subsets of $X$, and let $C(X)$ be its subspace of subcontinua. The Hausdorff metric on
$2^X$ is
\[
d_H(A,B)=\max\Bigl\{\sup_{a\in A}d(a,B),\ \sup_{b\in B}d(b,A)\Bigr\},
\qquad d(x,A)=\inf_{a\in A}d(x,a).
\]
We recall the following standard hyperspace facts.

\begin{lemma}[{\cite[Theorems 4.13 and 4.17]{ref14}}]
The space $(2^X,d_H)$ is compact, and $C(X)$ is closed in $2^X$.
\end{lemma}

A family of compact sets is a \emph{null family} if, for every $\varepsilon>0$, only
finitely many of its members have diameter at least $\varepsilon$.

\begin{lemma}
For every $\varepsilon>0$ there is a finite bound $P_X(\varepsilon)$ on the
cardinality of a family of pairwise disjoint connected subsets of a regular curve $X$
having diameter at least $\varepsilon$.
\end{lemma}

\begin{proof}
Choose a finite cover by open sets of diameter less than $\varepsilon$ with finite
boundaries. Every connected set of diameter at least $\varepsilon$ meets the union of
these boundaries, a finite set. We take $P_X(\varepsilon)$ to be the least possible
uniform bound.
\end{proof}

\begin{theorem}[{\cite[Lemma 5.6]{ref18}}]
Every nondegenerate Suslinian metric continuum has dimension one.
\end{theorem}

\begin{theorem}[{\cite[Theorem 1.6]{ref12}}]
A one-dimensional metric continuum $C$ is regular if and only if, for each
$\varepsilon>0$, there is an integer $N(\varepsilon)$ such that every family of
pairwise disjoint subcontinua of diameter at least $\varepsilon$ has at most
$N(\varepsilon)$ members.
\end{theorem}

\begin{corollary}
If a nondegenerate metric continuum satisfies the uniform bound in Theorem 2.5, then
it is regular.
\end{corollary}

\begin{proof}
The bound makes the continuum Suslinian. Apply Theorems 2.4 and 2.5.
\end{proof}

The bound is uniform over all disjoint families. The weaker condition that every
disjoint sequence of subcontinua is a null family is called the finitely Suslinian
property.

\subsection{The component quotient.}
In this subsection $M$ is a minimal set of a $G$-action on a regular curve $X$. The
following proposition records the properties of the quotient by connected components.

\begin{proposition}
Let $M$ be a minimal set and let $q\colon M\to Q$ be the quotient by its connected
components. Then $Q$ is the maximal zero-dimensional factor of $M$. Either $Q$ is
finite, with $G$ transitive on its fibers and the stabilizer of each fiber acting
minimally on it, or $Q$ is a Cantor set. In the latter case $q$ is almost one-to-one,
its nondegenerate fibers form a countable null family, and every invariant probability
on $Q$ lifts uniquely to $M$. The lift and its projection are measure-theoretically
isomorphic via $q$.
\end{proposition}

\begin{proof}
In a compact Hausdorff space components coincide with intersections of their clopen
neighborhoods. A compact metric space has at most countably many clopen sets, since
each is a finite union of elements of a countable basis. Membership in these sets
identifies $Q$ with a compact subspace of a countable product of discrete two-point
spaces. Every map to a zero-dimensional compact space is constant on components. The
action on $Q$ is minimal. If $Q$ is finite, its fibers are clopen, and
$Gx\cap C=G_Cx$ for $x\in C$, where $G_C=\{g:gC=C\}$. If $Q$ is infinite, minimality
excludes isolated points, so $Q$ is a Cantor set. Lemma 2.3 bounds the number of large
fibers. The exceptional values form a countable set and its complement is therefore
dense. Every invariant probability on an infinite minimal space is nonatomic. Off the
exceptional values the inverse of $q$ is continuous in the relative topology by
compactness, which gives the unique lift and the measure isomorphism.
\end{proof}

\subsection{Structure results for regular-curve actions.}
The following results will be used in the classification of minimal sets admitting an
invariant probability measure.

\begin{theorem}[{\cite[Theorem 2.3]{ref10}}]
Every action on a regular curve is null and hence tame. In particular, the action on
every closed invariant subset is tame.
\end{theorem}

\begin{theorem}[{\cite[Corollary 5.4(2)]{ref9}}]
Let a group act minimally on a compact metric space $M$. If the action is tame and
admits an invariant probability measure, then the factor map from $M$ to its maximal
equicontinuous factor is almost one-to-one.
\end{theorem}

\begin{theorem}[{\cite[Theorem 2.12]{ref11}}]
For an action on a locally connected compact metric space, every fiber of the maximal
equicontinuous factor map is connected.
\end{theorem}

\begin{theorem}[{\cite[Theorem 3.6]{ref4}}]
Every homogeneous nondegenerate Suslinian continuum is homeomorphic to the circle.
\end{theorem}
\section{Invariant neighborhoods and finite orbits}

In this section $G$ acts on a regular curve $X$. We prove Theorem A. We first
construct finite invariant sets separating invariant compact sets. We then use these
sets to obtain finite permutation actions on the components approaching a minimal set.

\begin{lemma}
Let $A$ be a nonempty compact invariant subset of $X$, and let $M$ be an infinite
minimal set disjoint from $A$. There is an invariant open set $U$ such that
$A\subset U$, $U\cap M=\emptyset$, $0<|\partial U|<\infty$.

If $A$ is connected, $U$ may be chosen connected.
\end{lemma}

\begin{proof}
Call an open set $V$ admissible if $A\subset V$ and $V\cap M=\emptyset$. A finite union
of finite boundary neighborhoods covering $A$ gives an admissible set with finite
boundary. Let $m$ be the least boundary cardinality among such sets. As $X$ is
connected, $m\ge1$. For open $V,W$ with finite boundary,
\[
|\partial(V\cup W)|+|\partial(V\cap W)|\le|\partial V|+|\partial W|.
\]
Indeed, both boundaries on the left lie in $\partial V\cup\partial W$. A point
belonging to just one of $\partial V,\partial W$ can contribute to at most one term on
the left, since it has a neighborhood in the other open set or disjoint from its
closure. Both $V\cup W$ and $V\cap W$ are admissible if $V,W$ are admissible. Thus a
finite union of admissible sets with $m$ boundary points has exactly $m$ boundary
points. If $A$ is connected, the component containing $A$ of any admissible set with
$m$ boundary points is open by local connectedness and its boundary is contained in
that of the original set. It too is admissible and has $m$ boundary points. In this
case take only the connected admissible sets with $m$ boundary points. In the general
case take all admissible sets with $m$ boundary points. The union $U$ of the chosen
family is invariant. In the connected case it is connected since all members contain
$A$. We have $|\partial U|\le m$. Otherwise take pairwise disjoint connected open
neighborhoods of $m+1$ distinct points of $\partial U$. Each neighborhood meets a
member of the family, so they all meet some finite union $V$ of members. Each also
contains a point of $X\setminus U$. By connectedness they all meet $\partial V$,
contrary to $|\partial V|=m$. Finally, $U\cap M=\emptyset$ and $\partial U$ is finite
and invariant. Were $\partial U$ to meet $M$, minimality would force the infinite set
$M$ into $\partial U$. Hence $U\cap M=\emptyset$. The boundary is nonempty by
connectedness of $X$.
\end{proof}

\begin{theorem}
Let $A$ be a nonempty compact connected invariant subset of a regular curve $X$. Every
infinite minimal set $M$ disjoint from $A$ is a Cantor set on which the action is
equicontinuous.
\end{theorem}

\begin{proof}
Let $S$ be the union of all connected invariant open neighborhoods of $A$ whose
closures miss $M$. Lemma 3.1 shows that the family is nonempty. The union is
invariant, open and connected, and is disjoint from $M$. If $S$ missed $M$, the same
lemma applied to $S$ would give a member of the family containing $S$. This would make
$S$ both open and closed in $X$, a contradiction. Thus $\partial S$ meets $M$ and, by
invariance and minimality, $M\subset\partial S$.

\smallskip
\noindent\emph{Claim 1.} The set $S$ admits an exhaustion by compact connected
invariant sets whose interiors are nested. Choose a countable subfamily
$(V_n)_{n\ge1}$ covering $S$. There exist compact connected invariant sets $K_n$ such
that
\[
K_n\subset\operatorname{int}K_{n+1}\subset K_{n+1}\subset S,
\qquad S=\bigcup_{n\ge1}K_n .
\]
Indeed, the closure of each member of the family lies in $S$, by the lemma applied to
that closure. Start with $K_1=\overline{V_1}$. Given $K_n$, apply the lemma to
$K_n\cup V_{n+1}$, which is connected since both sets contain $A$. The closure of the
resulting connected invariant neighborhood gives $K_{n+1}$.

\smallskip
\noindent\emph{Claim 2.} At each level only finitely many complementary components
reach $\partial S$, and the mesh of the complementary components tends to zero. Let
$\mathcal{V}_n$ be the collection of components of $S\setminus K_n$ whose closures
meet $\partial S$. Components of $S\setminus K_n$ are open and arcwise connected. Each
has closure meeting $K_n$, since otherwise it would be both open and closed in $S$.
Write $d_n=d(K_n,\partial S)>0$. A component whose closure meets $\partial S$ has
diameter at least $d_n$. More generally, there are only finitely many components of
$S\setminus K_n$ of diameter at least $d_n/2$, since each contains a compact arc of
diameter at least $d_n/3$ and these arcs are disjoint. This follows from Lemma 2.3. In
particular $\mathcal{V}_n$ is finite. For any $x\in\partial S$, take a sequence in
$S$ converging to $x$. Eventually its terms lie in components of $S\setminus K_n$ of
diameter at least $d_n/2$. Thus $x$ belongs to the closure of some member of
$\mathcal{V}_n$. We also have
\[
\sup\{\diam V: V \text{ is a component of } S\setminus K_n\}\longrightarrow 0 .
\]
Otherwise, for arbitrarily large $n$ there would be a compact arc of diameter at least
a fixed $\varepsilon>0$ in $S\setminus K_n$. Since every compact subset of $S$ is
contained in some $K_j$, these arcs could be chosen successively disjoint, contrary to
Lemma 2.3.

\smallskip
\noindent\emph{Claim 3.} The frontier $\partial S$ is a factor of an inverse limit of
finite $G$-sets. The map taking $V\in\mathcal{V}_{n+1}$ to its containing component in
$\mathcal{V}_n$ is equivariant. Let $E$ be the inverse limit of these finite $G$-sets.
For each $(V_n)\in E$ the nested sets $\overline{V}_n$ have exactly one common point.
It lies in $\partial S$, since the interiors of $K_n$ exhaust $S$. This defines a
continuous equivariant map $E\to\partial S$, continuous by the uniform diameter bound.
It is onto. For each $x\in\partial S$ and every $N$, choose $V_N\in\mathcal{V}_N$ with
$x\in\overline{V}_N$ and take its containing components at the earlier levels.
Compactness of the product of the finite sets gives an infinite compatible sequence
with $x$ in every closure.

Let $L$ be the closure of the image of $G$ in $\prod_n\Sym(\mathcal{V}_n)$. This
profinite group acts continuously on $E$. The closed fiber relation of $E\to\partial S$
is $L$-invariant by density of the image of $G$, so the action descends continuously
to $\partial S$. For $x\in M$, we have $M=\overline{Gx}=Lx$. Thus $M$ is equicontinuous
and homeomorphic to $L/L_x$, where $L_x$ is closed. This quotient is zero-dimensional:
the sets $NL_x$, for open normal subgroups $N$ of $L$, are clopen saturated
neighborhoods forming a basis at $L_x$. Translates give such a basis everywhere. As
$M$ is infinite and minimal, it has no isolated points and is a Cantor set.
\end{proof}

\begin{corollary}
If a group acting on a regular curve has a finite orbit, each infinite minimal set is a
Cantor set with an equicontinuous action. If the action has two distinct minimal sets,
then it has a finite orbit.
\end{corollary}

\begin{proof}
Let $H$ be the finite-index normal subgroup fixing every point of a given finite orbit,
and let $M$ be an infinite $G$-minimal set. If $N\subset M$ is $H$-minimal, its
finitely many distinct translates are $H$-minimal and are pairwise disjoint. They form
a clopen partition of $M$. Each translate is infinite, or their union $M$ would be
finite. The group $H$ has a fixed point outside $M$, so Theorem 3.2 makes its action on
each partition member equicontinuous and each member a Cantor set. This gives
equicontinuity of $H$ on $M$, since there are only finitely many clopen members. Taking
finitely many coset representatives gives equicontinuity of $G$ as well, and $M$ is a
Cantor set. If either of two distinct minimal sets is finite, the second assertion is
immediate. Otherwise apply Lemma 3.1 with one minimal set as $A$ and the other as $M$.
The finite nonempty invariant boundary contains a finite orbit.
\end{proof}

\begin{theorem}
Let $M$ be an infinite minimal set for a countable group $G$ acting on a regular curve
$X$. If $M$ has a nontrivial finite factor, then $G$ has a finite orbit in $X$ and its
action on $M$ is equicontinuous. In particular, if the maximal equicontinuous factor
of a Cantor minimal set is a Cantor set, then the factor map is a conjugacy.
\end{theorem}

\begin{proof}
Let $H$ be the kernel of the permutation action on the finite factor. Two distinct
factor fibers are disjoint nonempty compact $H$-invariant sets. Each contains an
$H$-minimal set. By Corollary 3.3, $H$ has a finite orbit. Since $H$ has finite index,
$G$ has a finite orbit as well. Corollary 3.3 gives equicontinuity on $M$.

To obtain the last assertion, let $Z$ be the Cantor maximal equicontinuous factor and
choose a compatible invariant metric on $Z$. Take a proper nonempty clopen subset $V$
of $Z$. Its distance from $Z\setminus V$ is positive and is the same for every
translate. Consequently the translates of $V$ are determined by their traces on a fixed
sufficiently fine finite net. They generate a finite $G$-invariant clopen partition
with at least two atoms. This gives a nontrivial finite factor of $M$. As $M$ is
equicontinuous, its maximal equicontinuous factor map is injective.
\end{proof}

\begin{corollary}
If the ambient action has no finite orbit, it has the same unique minimal set under
every finite-index subgroup. An infinite minimal set with finitely many components is
connected.
\end{corollary}

\begin{proof}
Let $H\le G$ have finite index and take a normal subgroup $N$ of finite index contained
in $H$. A finite $N$-orbit would give a finite $G$-orbit. Thus $N$ has a unique minimal
set $L$ by Corollary 3.3. Normality makes $L$ invariant under $G$. It is minimal under
both $G$ and $H$ and unique for both, since every nonempty compact invariant set for
either contains an $N$-minimal set. Finally, if an infinite minimal set has a finite
number greater than one of components, its component quotient is a nontrivial finite
factor. By Theorem 3.4 and Corollary 3.3 it is a Cantor set, which has infinitely many
components. This contradiction proves the assertion.
\end{proof}

In particular, an action with a global fixed point cannot have a disjoint nondegenerate
connected minimal set. No invariant measure assumption on that set is needed.

\begin{proof}[Proof of Theorem A]
The finite-orbit assertion and uniqueness follow from Corollary 3.3. The assertions
concerning finite-index subgroups and finitely many components are given by Corollary
3.5.
\end{proof}
\section{Minimal sets with an invariant measure}

In this section we prove the classification part of Theorem B. We first study the
fibers of almost one-to-one equicontinuous factors, and then apply the structure
results recalled in Section 2. The realization of the four possibilities will be given
in Section 5.

\subsection{Fibers of almost automorphic minimal sets.}
Let a group $G$ act by homeomorphisms on a regular curve $(X,d)$, and let $M\subset X$
be a minimal set. In this subsection, $\pi\colon M\to Z$ is a factor map onto a minimal
equicontinuous system. Write $K_z=\pi^{-1}(z)$ and suppose that $\pi$ is almost
one-to-one. Thus the set $D=\{z\in Z:|K_z|=1\}$ is dense in $Z$, and
$M_0=\pi^{-1}(D)$ is dense in $M$. Indeed, $M_0$ is nonempty and invariant. The set $D$
is a $G_\delta$ since $z\mapsto\diam K_z$ is upper semicontinuous.

\begin{theorem}
Every fiber of $\pi$ is finite. For each $\varepsilon>0$,
\[
\bigl|\{z\in Z:\diam K_z\ge\varepsilon\}\bigr|\le P_X(\varepsilon).
\]
In particular, the nonsingleton fibers form a null family, and there are at most
countably many of them. If some point $b\in M$ has order at most $n\ge1$ in $X$, then
$|K_z|\le n$ for every $z\in Z$.
\end{theorem}

Equip $Z$ with a compatible invariant metric and let $L$ be the compact closure of the
image $\rho(G)$ in its isometry group. The action of $L$ on $Z$ is transitive. Each
orbit map $L\to Z$ is open, since it identifies $Z$ with the quotient by a closed
stabilizer. If $\rho(u_j)\to e$, then $u_jp\to p$ for every $p\in M_0$, by compactness
and the uniqueness of the lift of $\pi(p)$.

\begin{lemma}
For $z\in Z$ and $b\in K_z$, there exist $u_j\in G$ such that $\rho(u_j)\to e$ and
$u_jK_z\to\{b\}$ in the Hausdorff metric.
\end{lemma}

\begin{proof}
Choose $p\in M_0$ and $h_j\in G$ with $h_jp\to b$. After passing to a subsequence,
suppose $\rho(h_j)\to h$. Then $h\pi(p)=z$ and $\rho(h_j)^{-1}z\to\pi(p)$. Upper
semicontinuity gives $h_j^{-1}K_z\to\{p\}$. For each fixed $j$ choose $m(j)>j$ so large
that $d_H(h_jh_{m(j)}^{-1}K_z,\{h_jp\})<1/j$. Set $u_j=h_jh_{m(j)}^{-1}$.
\end{proof}

\begin{proof}[Proof of Theorem 4.1]
\emph{Claim 1.} Every fiber is finite. Suppose some fiber $K_z$ is infinite. Let $a$ be
an accumulation point of $K_z$ and choose $b\in K_z\setminus\{a\}$. Fix
$0<\delta<d(a,b)/5$. For arbitrary $r$, choose distinct $x_1,\dots,x_r\in K_z\cap
B(a,\delta)$. By local connectedness there are pairwise disjoint connected open
neighborhoods $U_i\subset B(a,2\delta)$ of $x_i$. Choose $p_i\in M_0\cap U_i$. For the
sequence in Lemma 4.2, we have $u_jx_i\to b$ and $u_jp_i\to p_i$. Since
$d(p_i,b)>3\delta$, for all sufficiently large $j$ the disjoint connected sets $u_jU_i$
each have diameter greater than $\delta$. This gives $r\le P_X(\delta)$, a
contradiction.

\smallskip
\noindent\emph{Claim 2.} At most $P_X(\varepsilon)$ fibers have diameter at least
$\varepsilon$. For distinct $z_1,\dots,z_r$ with $\diam K_{z_i}\ge\varepsilon$, choose
$x_i,y_i\in K_{z_i}$ at distance at least $\varepsilon$. Each set
$\{h\in L:hz_i\in D\}$ is dense $G_\delta$ by openness of the orbit map. Choose $h$ in
their intersection, and choose $g_j\in G$ with $\rho(g_j)\to h$. For each $i$ both
$g_jx_i$ and $g_jy_i$ converge to the unique lift of $hz_i$. These $r$ limiting points
are distinct. Choose disjoint connected open neighborhoods of them in $X$. For large
$j$ each pair is in its corresponding neighborhood. Pulling the neighborhoods back by
$g_j$ gives $r$ disjoint connected subsets of $X$, each of diameter at least
$\varepsilon$. Hence $r\le P_X(\varepsilon)$. Taking $\varepsilon=1/k$ proves
countability.

\smallskip
\noindent\emph{Claim 3.} If $b\in M$ has order at most $n$ in $X$, then every fiber has
at most $n$ points. Let $b\in M$ have order at most $n\ge1$, and fix a fiber $K_z$ with
$|K_z|>1$. Choose $x\in K_z$. By minimality there are $h_j\in G$ with $h_jx\to b$.
Lemma 4.2 gives $v_mK_z\to\{x\}$. Choosing $m(j)$ sufficiently large for each $j$, we
obtain $g_j=h_jv_{m(j)}$ with $g_jK_z\to\{b\}$. In particular,
$\rho(g_j)z\to\pi(b)$.

For distinct $x_1,\dots,x_r\in K_z$, take disjoint connected open neighborhoods $U_i$
of $x_i$ and $p_i\in M_0\cap U_i$. Since $K_z$ is not a singleton, $\pi(p_i)\neq z$
for each $i$. Let $d_Z$ be the invariant metric on $Z$ and set
$\eta=\tfrac13\min_i d_Z(\pi(p_i),z)>0$. For every $j$ and $i$,
\[
d_Z(\pi(g_jp_i),\rho(g_j)z)=d_Z(\pi(p_i),z)\ge3\eta.
\]
Consequently $\pi(g_jp_i)\notin B(\pi(b),2\eta)$ for all large $j$. Choose an open
neighborhood $V$ of $b$ with $|\partial_X V|\le n$ and
$V\cap M\subset\pi^{-1}(B(\pi(b),\eta))$. Then $g_jx_i\in V$ and $g_jp_i\notin V$ for
every $i$ and all sufficiently large $j$. The disjoint connected sets $g_jU_i$ all meet
$\partial_X V$, and $r\le n$.
\end{proof}

\subsection{Connected minimal sets and the measured classification.}
We first consider a connected minimal set supporting an invariant probability measure.

\begin{lemma}
Let $H$ act minimally on a nondegenerate regular curve $C$. If $C$ supports an
$H$-invariant probability measure, then $C$ is a circle and the action is conjugate to
a group of circle isometries.
\end{lemma}

\begin{proof}
By Theorems 2.8 and 2.9, the maximal equicontinuous factor map $\pi\colon C\to Z$ is
almost one-to-one. Its fibers are finite by Theorem 4.1. On the other hand, Theorem
2.10 shows that they are connected. Hence $\pi$ is a homeomorphism. The closure of the
acting group for an invariant compatible metric is a transitive compact group of
isometries. Thus $C$ is homogeneous. It is Suslinian by Lemma 2.3. It follows from
Theorem 2.11 that $C$ is a circle. The invariant probability has full support and is
nonatomic by minimality. Measuring oriented arcs from a fixed origin gives a circle
coordinate in which this measure is Lebesgue measure. Homeomorphisms preserving
Lebesgue measure on the circle are rotations or reflections.
\end{proof}

\begin{theorem}
Let $M$ be a minimal set of an action on a regular curve $X$. If $M$ supports an
invariant probability measure, then it is a finite orbit, a Cantor set or a circle with
an isometric action in a suitable coordinate. In the circle case $M$ is the unique
minimal set of $X$.
\end{theorem}

\begin{proof}
Theorems 2.8 and 2.9 show that the maximal equicontinuous factor map $\pi\colon M\to Z$
is almost one-to-one. Theorem 4.1 makes every fiber finite. Suppose that $C$ is a
nondegenerate component with infinitely many distinct translates $g_jC$. Lemma 2.3
gives $\diam(g_jC)\to0$. Since $Z$ admits a compatible invariant metric,
\[
\diam_Z\pi(C)=\diam_Z\pi(g_jC)\longrightarrow0 .
\]
This would put $C$ in one finite fiber. Thus $C$ has only finitely many translates,
whose union is all of $M$. Corollary 3.5 gives $M=C$. Lemma 4.3 gives the circle
isometry action. Corollary 3.3 excludes any second minimal set, as its existence would
make $M$ a Cantor set. If $M$ has only singleton components, Proposition 2.7 gives the
other two possibilities.
\end{proof}

\subsection{The maximal equicontinuous factor of a Cantor minimal set.}
Suppose $M$ is a Cantor minimal set in $X$ supporting a $G$-invariant probability
measure. Let $\pi\colon M\to Z$ be its maximal equicontinuous factor map. By Theorem
4.1, its fibers form a null family, in the sense that only finitely many have diameter
at least a given positive number.

\begin{proposition}
There are a regular curve $Y$ with a $G$-action and a factor map $P\colon X\to Y$ whose
restriction to $M$ is $\pi$, after identifying $Z$ with a closed invariant subset of
$Y$. The restriction of $P$ to $X\setminus M$ is a homeomorphism onto $Y\setminus Z$.
\end{proposition}

\begin{proof}
On $X$, take the equivalence relation consisting of the diagonal and all pairs
$(u,v)\in M^2$ satisfying $\pi(u)=\pi(v)$. This is a closed invariant relation, so its
quotient $Y$ is a compact metric continuum with an induced $G$-action. The quotient map
$P$ is closed, and its singleton fibers outside $M$ give the claimed homeomorphism. The
restriction of $P$ to $M$ identifies its image with $Z$. To check regularity, let
$y\in Y$ and let $V$ be an open neighborhood of $y$. If $y\notin Z$, take a
finite-boundary open neighborhood of its preimage whose closure is contained in
$P^{-1}(V)\setminus M$. Its image is an open neighborhood of $y$ with finite boundary.
For $y\in Z$, choose a clopen subset $A$ of $M$ with
$P^{-1}(y)\subset A\subset P^{-1}(V)\cap M$. Cover $A$ by finitely many open subsets of
$X$ with finite boundary and with closures contained in
$P^{-1}(V)\setminus(M\setminus A)$. Their union $U$ satisfies $U\cap M=A$ and
$\partial_X U\cap M=\emptyset$. The set
\[
B_A=\pi(A)\cap\pi(M\setminus A)
\]
is finite. Indeed, if $A$ and $M\setminus A$ are both nonempty, the fiber over every
point of $B_A$ has diameter at least $d(A,M\setminus A)>0$. Now $W=Y\setminus P(X\setminus U)$
is an open neighborhood of $y$ contained in $V$, and
$\partial_Y W\subset P(U)\cap P(X\setminus U)\subset P(\partial_X U)\cup B_A$. This
proves regularity.
\end{proof}

\begin{proposition}
The system $Z$ is either a Cantor minimal equicontinuous system or a minimal circle
isometry action. In the Cantor case $\pi$ is a conjugacy and $G$ has a finite orbit in
$X$. In the circle case, all fibers of $\pi$ have at most two points. The values with
two-point fibers form a countable dense invariant set $D$. In that case $M$ is the
unique minimal set in $X$. The space $M$ is obtained from the circle by replacing each
$t\in D$ by two adjacent points in the circular order. Each element of $G$ acts by the
unique lift of its circle isometry.
\end{proposition}

\begin{proof}
The system $Z$ is infinite: by the almost one-to-one property, a finite minimal $Z$
would consist entirely of singleton fibers. Apply Theorem 4.4 to the minimal set $Z$ in
the regular curve $Y$ of Proposition 4.5. The required invariant probability is the
image of the one on $M$. This gives the two possibilities for $Z$. Theorem 3.4 settles
the Cantor case.

Suppose $Z$ is a circle. By Theorem 4.4 it is the unique minimal set in $Y$. Every
compact invariant subset of $X$ disjoint from $M$ would have image disjoint from $Z$.
This proves uniqueness of $M$. Fix $t\in Z$ and choose a circle coordinate near $t$.
For a clopen $A\subset M$, the two closed sets $\pi(A)$ and $\pi(M\setminus A)$ cover
$Z$ and meet only in the finite set $B_A$ of the preceding proof. A sufficiently short
open interval to the right of $t$ avoids $B_A$. By connectedness, all its fibers are
contained in $A$, or all are contained in $M\setminus A$. It follows that lifts of
points approaching $t$ from the right have only one possible limit in $M$. Otherwise, a
clopen set separating two possible limits would contradict this property. Compactness
gives existence of the limit, which we denote by $t^+$. The left-hand limit $t^-$
exists in the same sense. Density of the union of singleton fibers gives
$\pi^{-1}(t)=\{t^-,t^+\}$. Indeed, each point of a nonsingleton fiber is approximated
by singleton lifts outside that fiber; their images admit a subsequence approaching $t$
from one side. The assertion also holds for singleton fibers by continuity of $\pi$.

The set $D=\{t:t^-\neq t^+\}$ is countable by Theorem 4.1. It is nonempty since $M$ is
a Cantor set, and is dense by invariance and minimality. In the circle with points of
$D$ split, points with base coordinates approaching $t$ from the left converge to the
left copy of $t$, and points approaching from the right converge to the right copy. In
the circular order topology these one-sided limits and convergence of eventually
constant sequences in a fiber determine convergence. They give exactly the limits just
proved in $M$. Thus matching the left and right copies to $t^-$ and $t^+$ gives a
homeomorphism. Here the split circle is metrizable, as can also be seen by inserting
disjoint intervals of summable positive lengths at the countably many split values and
retaining their endpoints. A lift of a circle isometry is determined on singleton
fibers, hence on their dense union and on $M$.
\end{proof}

\subsection{Cyclic subactions.}
Throughout this subsection, $G$ acts on a regular curve $X$, and $M$ is a $G$-minimal
set supporting a $G$-invariant probability measure $\mu$. For $f\in G$, we write
$\Per_X(f)$ for the set of its periodic points in $X$. A point is nonwandering for $f$
if each of its neighborhoods meets a positive iterate of itself. Since $\mu$ has full
support on $M$, every point of $M$ is nonwandering for $f$. Indeed, each relatively
open nonempty subset of $M$ has positive measure and meets a positive iterate of itself.

We recall the results for individual homeomorphisms which are needed below. For a
homeomorphism $f$, write $\omega_f(x)$ for the set of limits of subsequences $f^{n_j}x$
with $n_j\to\infty$, and write $\Omega(f)$ for its nonwandering set. A point $x$ is
recurrent if $x\in\omega_f(x)$. An adding machine is a minimal translation on an
infinite inverse limit of finite cyclic groups.

\begin{theorem}[{\cite[Theorems 1.1(i), 2.1 and 4.1]{ref15}}]
Let $f$ be a homeomorphism of a regular curve. Then every $\omega$-limit set is
minimal, and every nonwandering point is recurrent. If $f$ has no periodic point, then
it has a unique minimal set.
\end{theorem}

\begin{theorem}[{\cite[Theorem 4.1]{ref2}}]
Let $f$ be a homeomorphism of a regular curve $X$ with no periodic point. There are an
irrational number $\alpha$ and a continuous surjection $P\colon X\to S^1$ with
connected fibers such that $P\circ f=R_\alpha\circ P$, $R_\alpha(t)=t+\alpha\pmod 1$.
\end{theorem}

\begin{theorem}[\cite{ref3}]
If a homeomorphism of a regular curve has a periodic point, then its restriction to
every infinite minimal set is conjugate to an adding machine.
\end{theorem}

\begin{lemma}[{\cite[Lemma 2]{ref5}}]
Let $f$ be a homeomorphism of a compact zero-dimensional metric space $K$. If $K$ is
the union of periodic orbits and invariant subsets on which $f$ is conjugate to an
adding machine, then $\{f^n:n\in\mathbb{Z}\}$ is equicontinuous.
\end{lemma}

\begin{theorem}
Suppose some $f\in G$ satisfies $\Per_X(f)=\emptyset$. Then $M$ is the unique minimal
set of $f$. There are a homomorphism $\rho\colon G\to\Isom(S^1)$ and a factor map
$p\colon M\to S^1$ such that $\rho(f)$ is an irrational rotation. The map $p$ is the
maximal equicontinuous factor map both for $f$ and for $G$. Every fiber has one or two
points. The set of values with two-point fibers is at most countable, and the singleton
fibers have dense union in $M$. Moreover, $\ker\rho$ consists precisely of the elements
acting trivially on $M$. If $p$ is injective, then $M$ is a circle. Otherwise $M$ is a
Cantor set, and the two-point fibers lie over a dense invariant subset of the circle.
\end{theorem}

\begin{proof}
Every point of $M$ is nonwandering for $f$. By Theorem 4.7, it belongs to the unique
$f$-minimal set. This gives the first assertion. Theorem 4.8 gives a monotone surjection
$P\colon X\to S^1$ such that $Pf=R_\alpha P$ for an irrational rotation $R_\alpha$. Put
$p=P|_M$, which is onto by minimality of $R_\alpha$. The nondegenerate fibers of $P$
are disjoint subcontinua, so they are at most countably many and form a null family.
The lifts of all singleton fibers of $P$ belong to $M$, and their union is dense in $M$
by $f$-minimality. If $p(x)=p(x')=t$, then $f^nx$ and $f^nx'$ belong to the same fiber
$P^{-1}(R_\alpha^n t)$. These fibers are pairwise disjoint for distinct $n$, so
$d(f^nx,f^nx')\to0$ as $n\to\infty$. Consequently every equicontinuous factor of
$(M,f)$ is constant on the fibers of $p$. Thus $p$ gives its maximal equicontinuous
factor.

Let $\pi\colon M\to Z$ be the maximal equicontinuous factor map for $G$. By Theorems
2.8 and 2.9, the map $\pi$ is almost one-to-one. Since $Z$ is equicontinuous also under
$f$, there is a factor map $r$ for the cyclic actions such that $\pi=rp$. The
equivalence relation defined by $r$ is closed and invariant under simultaneous
translation of both coordinates by $\alpha$. Density of the powers of $R_\alpha$ makes
it invariant under every circle translation. All fibers of $r$ are therefore translates
of one another. As $\pi$ has a singleton fiber, $r$ has one as well, and $r$ is a
homeomorphism. Transferring the $G$-action on $Z$ to the circle makes $p$ a $G$-factor
map. The measure $p_*\mu$ is invariant under $R_\alpha$, hence under all rotations by
density and continuity, and is therefore Lebesgue measure. The induced circle
homeomorphisms preserve this measure and are therefore isometries. Any element inducing
the identity on the circle fixes the dense union of singleton fibers pointwise. By
Theorem 4.4, $M$ is a circle or a Cantor set. In the former case the $G$-action is
equicontinuous, and its maximal equicontinuous factor map is injective. In the latter
case Proposition 4.6 gives precisely the asserted description of its fibers.
\end{proof}

\begin{corollary}
Suppose $M$ is a Cantor set and its maximal equicontinuous factor is a circle. For
every $f\in G$, the set $\Per_X(f)$ is empty if and only if $f$ induces an irrational
rotation on the circle. If $G$ is finitely generated and $M$ is a Cantor set, its
maximal equicontinuous factor is a circle if and only if $G$ contains an element with
no periodic point in $X$.
\end{corollary}

\begin{proof}
Write $p\colon M\to S^1$ for the factor map. Its fibers are singletons on a dense
subset of $M$. If the isometry induced by $f$ has finite order, a positive power of $f$
fixes all singleton lifts, hence fixes $M$ pointwise. Suppose instead that $f$ induces
an irrational rotation. Every $f$-orbit closure in $M$ projects onto the circle and
therefore contains the dense union of singleton fibers. Thus $f$ is minimal on $M$. By
Proposition 4.6, there is a two-point fiber. Its points are asymptotic under $f$, since
its translates are distinct fibers of a null family. Hence $(M,f)$ is not an adding
machine. Theorem 4.9 therefore excludes an ambient periodic point. Finally, if $G$ is
finitely generated, so is the subgroup of index at most two preserving orientation on
the circle factor. Its rotation image is infinite by minimality and therefore contains
an element of infinite order. The converse follows from Theorem 4.11.
\end{proof}

\begin{proposition}
Suppose $M$ is a Cantor set. If $f\in G$ has a periodic point in $X$, then the family
$\{f^n|_M:n\in\mathbb{Z}\}$ is equicontinuous.
\end{proposition}

\begin{proof}
Since every point of $M$ is nonwandering for $f$, Theorems 4.7 and 4.9 imply that each
point belongs to a periodic orbit or an adding machine. Lemma 4.10 now gives
equicontinuity.
\end{proof}

\begin{corollary}
Suppose $G$ is finitely generated and $M$ is a Cantor minimal set supporting an
invariant probability. Then the action on $M$ is equicontinuous if and only if every
element of $G$ has a periodic point in $X$.
\end{corollary}

\begin{proof}
By Proposition 4.6 the equicontinuous case is exactly the one with a Cantor maximal
equicontinuous factor. Corollary 4.12 now gives the equivalence. The rational rotation
example in Subsection~\ref{subsec:split} shows why finite generation is needed for this
criterion.
\end{proof}
\section{Models and examples}

Each alternative in Theorem B occurs. A finite transitive $G$-set $E$ can be realized
by taking one unit interval for each point of $E$, identifying their initial endpoints,
and permuting the intervals by $G$. Their other endpoints form a minimal set conjugate
to $E$. A minimal group of circle isometries acts on the circle itself. We next give
the Cantor examples and examples without an invariant probability.

\subsection{Equicontinuous Cantor actions on dendrites.}

\begin{proposition}
Every minimal equicontinuous action of a countable group $G$ on a Cantor set $Z$ is
conjugate to the action on the set of endpoints of a dendrite. The dendrite admits a
compatible metric for which $G$ acts by isometries and fixes a point.
\end{proposition}

\begin{proof}
Take a compatible $G$-invariant metric on $Z$. Construct finite $G$-invariant clopen
partitions $\mathcal{P}_n$, starting with $\mathcal{P}_0=\{Z\}$, so that each atom
splits at the next level and the mesh tends to zero. Indeed, take a finite clopen
refinement of arbitrarily small mesh that splits every atom of the preceding partition.
Each clopen set has finitely many distinct translates: its positive distance from its
complement is preserved, so its translates are determined by their traces on a
sufficiently fine finite net. Taking the common refinement of all translates gives the
desired invariant partition.

Make each atom a vertex and join it to its parent atom at the preceding level. Give
edges from level $n-1$ to level $n$ length $2^{-n}$. Let $K_n$ be the finite subtree
through level $n$. Every descendant has distance at most $2^{-n}$ from $K_n$, so the
completion $X$ in the path metric is a compact real tree. Its new points are precisely
the infinite rays. Every vertex has at least two children, so exactly the ray endpoints
are endpoints of $X$. An infinite nested sequence of atoms has singleton intersection,
and two rays with the same first $n$ levels have distance at most $2^{1-n}$. This
identifies the endpoint set homeomorphically with $Z$. The permutations of the atoms
induced by $G$ preserve inclusion and extend linearly on edges to isometries of $X$
fixing the root. The identification of endpoints is equivariant.
\end{proof}

The same binary tree gives two different equicontinuous actions. For
$Z=\mathbb{Z}_2$, the compact group of $2$-adic integers, partition $Z$ by residues
modulo $2^n$. The translation $x\mapsto x+1$ acts on level $n$ as a cycle of length
$2^n$, so its action on the endpoints is the dyadic adding machine. Alternatively,
identify an end with its digits $(b_0,b_1,\dots)\in\prod_{j=0}^{\infty}\{0,1\}$,
starting with the least significant digit. The group $\bigoplus_{j\ge0}\mathbb{Z}/2\mathbb{Z}$
acts by flipping the digits in any prescribed finite set of positions. Each orbit meets
every cylinder, hence is dense, although every nonidentity acting element has order
two. For example, modulo $4$ the adding machine permutes the residues by $(0\,1\,2\,3)$,
whereas flipping $b_0$ gives $(0\,1)(2\,3)$.

\subsection{Splitting orbits on the circle.}\label{subsec:split}

Let $\rho\colon G\to\Isom(S^1)$ be a minimal action, and let $D\subset S^1$ be a
nonempty countable invariant set. We construct a Cantor minimal set whose double fibers
over the circle are exactly those above $D$. The construction also allows reflections.
The set $D$ is dense by minimality. Assign a length $\ell_t>0$ to each $t\in D$, with
$\sum_{t\in D}\ell_t<\infty$, and set $\ell_t=0$ elsewhere. Choose an origin outside
$D$ and identify the base circle with $[0,1]/(0\sim1)$ in its positive orientation. Put
\[
L=1+\sum_{t\in D}\ell_t,\qquad
a(t)=t+\sum_{\substack{s\in D,\ s<t}}\ell_s,\qquad
b(t)=a(t)+\ell_t .
\]
The ordered intervals $I_t=[a(t),b(t)]$ partition $S_b=[0,L]/(0\sim L)$, with the
interval at the origin taken only once. Indeed, $a(t)$ and $b(t)$ are the left and right
limits of the cumulative length at $t$. They define a continuous monotone surjection
$P\colon S_b\to S^1$ by $P(I_t)=\{t\}$. For $g\in G$, let $g_b$ send $I_t$ affinely onto
$I_{\rho(g)t}$, preserving or reversing orientation according as $\rho(g)$ does. On
singleton fibers this determines the map uniquely. These maps preserve or reverse
circular order and are bijections, hence are homeomorphisms. Composition agrees with the
group law since the affine map between two oriented intervals is determined by the
images of their endpoints.

Put $M=S_b\setminus\bigcup_{t\in D}(a(t),b(t))$ and $p=P|_M$. Every endpoint is a limit
of singleton lifts on its exterior side. Thus $M$ is the closure of the singleton lifts
and is perfect. The density of $D$ implies that $M$ contains no nondegenerate arc, so
$M$ is a Cantor set. Every nonempty closed invariant subset of $S_b$ projects onto the
base by minimality and contains all singleton lifts. Consequently $M$ is the unique
minimal set, and $D$ is exactly the set of its double-fiber values. Lebesgue measure
gives full measure to the complement of $D$, and lifting it by its unique inverse there
gives an invariant probability on $M$. The factor $p$ is maximal equicontinuous. For
$t\in D$ choose $g_n$ with $\rho(g_n)t$ pairwise distinct. Then
$\ell_{\rho(g_n)t}\to0$, so $a(t)$ and $b(t)$ form a proximal pair. Every equicontinuous
factor identifies them and hence factors through $p$. The lifted action on $M$ has
exactly the kernel of $\rho$. Indeed, an element acting trivially on $M$ acts trivially
on its factor, and the prescribed lift of the identity is the identity on $S_b$.

For an explicit Denjoy example, take $G=\mathbb{Z}$, $\alpha=(\sqrt5-1)/2$, and
$t_k=k\alpha\pmod1$ for $k\in\mathbb{Z}$ in the standard base coordinate. Split the
orbit $D=\{t_k:k\in\mathbb{Z}\}$ with
\[
\ell_{t_k}=\tfrac13\,2^{-|k|}.
\]
The total inserted length is one. One may cut the base at $1/2$ in the construction
above. The generator carries $I_{t_k}$ to $I_{t_{k+1}}$. Now let $G=\mathbb{Q}/\mathbb{Z}$
act by all rational rotations and split any nonempty countable union of its orbits, for
example $D=\mathbb{Q}/\mathbb{Z}$ or $D=(\mathbb{Q}/\mathbb{Z})\cup(\alpha+\mathbb{Q}/\mathbb{Z})$.
Every lifted rational rotation has the same finite order as its base rotation, yet the
action on $M$ is minimal and fails to be equicontinuous because it has distinct proximal
endpoint pairs.

\subsection{The finitary alternating group.}
Let $G$ be the group of even permutations of $\mathbb{N}$ of finite support. This
countable group is the increasing union of the finite alternating groups $A_n$.
Averaging over these groups and taking a weak limit gives an invariant probability on
every nonempty compact invariant set. It has no nontrivial finite quotient: the
restriction of a homomorphism to a fixed finite group is trivial on $A_n$ once $n\ge5$
and $|A_n|$ exceeds the order of the image, by simplicity of $A_n$. It also has no
nontrivial homomorphism to the circle isometry group. The orientation character is
trivial by the finite quotient assertion, and the rotation group is abelian, whereas
$G$ is perfect. Every minimal set of a $G$-action on a regular curve is a singleton.
Indeed, a nontrivial finite orbit gives a nontrivial finite quotient. An infinite
equicontinuous Cantor minimal set has a nontrivial finite factor by the proof of
Theorem 3.4. In each of the remaining measured cases there is an infinite circle
isometry factor. All three possibilities for a nontrivial minimal set are therefore
excluded. A faithful action is obtained by joining intervals of lengths $2^{-n}$ at one
endpoint and permuting them by their indices, preserving normalized distance along each
interval. This compact star is a dendrite, its center is fixed, and the orbit closure of
each other point contains the center.

\subsection{Actions without an invariant probability.}
We give a Cantor example and a circle example. Take the Cayley tree of
$F_2=\langle a,b\rangle$. Give an edge from level $n$ to level $n+1$ length $2^{-n-1}$,
with levels measured from the identity. Its metric completion $X$ is a compact dendrite.
Indeed, outside the finite subtree of depth $n$, the components have diameters at most
$2^{1-n}$. Left translation extends to the completion: a word of length $m$ changes the
level of a vertex by at most $m$, so the edge lengths change by factors between $2^{-m}$
and $2^m$. The ends $M$ are infinite reduced words. They form a Cantor set, with a
clopen basis of cylinders $[u]$ specified by nonempty finite reduced words $u$. Given
nonempty cylinders $[u]$, $[w]$, choose a letter $c$ for which $ucw^{-1}$ is reduced. If
$\xi\notin[w]$, the whole of $w^{-1}$ cannot cancel against $\xi$, so
$(ucw^{-1})\xi\in[u]$. This proves both minimality and extreme proximality, that is,
every proper closed set can be sent into any nonempty open set. No invariant probability
$\mu$ exists on $M$. The equality $a(M\setminus[a^{-1}])=[a]$ would give
$\mu([a])+\mu([a^{-1}])=1$. The same equality for $b$ would assign total mass two to four
disjoint cylinders partitioning $M$.

For the circle example, let $\mathrm{PSL}_2(\mathbb{Z})$ act on the projective line
$\mathbb{R}\cup\{\infty\}$. Iterating $x\mapsto x+1$ sends every real point to $\infty$
in the limit. The orbit of $\infty$ contains every rational point, by B\'ezout's
identity. Thus the action is minimal. Any probability invariant under the translation
must be $\delta_\infty$, by dominated convergence, and $x\mapsto-1/x$ interchanges $0$
and $\infty$. Hence the circle supports no invariant probability for this action. Thus
both a Cantor set and a circle can support a minimal action without an invariant
probability.

\begin{proof}[Proof of Theorem B]
Theorem 4.4 gives the finite, Cantor and circle cases, including the conjugacy to circle
isometries and uniqueness in the circle case. For a Cantor minimal set, Proposition 4.6
applies to its maximal equicontinuous factor. It gives either an equicontinuous Cantor
action or the stated split-circle action. It also gives the finite ambient orbit in the
former case, and uniqueness, the two-point exceptional fibers and the lifted action in
the latter case. For the converse, the finite and circle models were given at the start
of this section. Proposition 5.1 realizes every minimal equicontinuous Cantor action for
the same group, and Subsection~\ref{subsec:split} realizes every split-circle action.
These systems carry invariant probabilities: the uniform measure in the finite case, an
invariant measure for a compact equicontinuous action in the second case, and Lebesgue
measure in the circle case. In the split-circle case, Lebesgue measure lifts uniquely
across the countable exceptional set, as shown in the construction.
\end{proof}
\section{Minimal sets with infinitely many components}

\subsection{Homeomorphism types.}
For a compact metric space $M$, let $q_M\colon M\to Q_M$ denote the quotient by its
connected components. Write $\mathcal{T}(M)$ for the set of homeomorphism types of its
nondegenerate components. The quotient $Q_M$ is compact, metrizable and
zero-dimensional. In particular, if a component $C$ is contained in an open set
$O\subset M$, then there is a clopen set $H$ with $C\subset H\subset O$. Indeed, the
open set $Q_M\setminus q_M(M\setminus O)$ contains $q_M(C)$, and we may take $H$ to be
the inverse image of a clopen neighborhood contained in this open set. Every clopen
subset of $M$ is a union of components.

\begin{lemma}
Suppose that the nondegenerate components of a compact metric space $(M,d)$ form a null
family. Equip $Q_M$ with a compatible metric $\rho$. Let $U\subset M$ be clopen and let
$C\subset U$ be a component with empty interior. For each $\delta>0$ there exist a
partition of $U\setminus C$ into nonempty clopen subsets $U_1,U_2,\dots$ of $M$ and
points $c_j\in C$ such that
\[
\diam_\rho q_M(U_j)<\delta\quad (j\ge1),\qquad
\lim_{j\to\infty}\sup_{u\in U_j}d(u,c_j)=0 .
\]
For every nonempty relatively open subset $O$ of $C$, there are infinitely many $j$
with $c_j\in O$.
\end{lemma}

\begin{proof}
Choose a decreasing sequence of clopen neighborhoods $H_m$ of $C$ such that $H_1=U$ and
$\bigcap_m H_m=C$. Compactness gives $\sup_{u\in H_m}d(u,C)\to0$. Put
$A_m=H_m\setminus H_{m+1}$ and let $s_m$ be the supremum of the diameters of the
components contained in $A_m$, with $s_m=0$ if $A_m$ is empty. Since the shells are
disjoint and there are only finitely many components of diameter at least
$\varepsilon$ for each $\varepsilon>0$, we have $s_m\to0$. Each component of a
nonempty $A_m$ has an open neighborhood in $A_m$ of diameter less than $s_m+1/m$. By the
clopen neighborhood property above and compactness, $A_m$ has a finite clopen cover by
sets of diameter less than $s_m+1/m$. Refine this cover to a disjoint clopen partition
by taking successive differences. Intersect its members with the inverse images of a
finite clopen partition of $Q_M$ of $\rho$-mesh less than $\delta$, and discard the
empty intersections. Enumerate the members of the resulting partitions in order of
increasing $m$. There are infinitely many, since otherwise $C$ would be open. For each
member $U_j$ choose $u_j\in U_j$ and $c_j\in C$ with $d(u_j,c_j)=d(u_j,C)$. If
$U_j\subset A_m$, then
\[
\sup_{u\in U_j}d(u,c_j)\le s_m+1/m+\sup_{v\in H_m}d(v,C).
\]
Each shell contributes only finitely many sets, so the right-hand side tends to zero as
$j\to\infty$.

Given $c\in C$, choose $x_n\in U\setminus C$ tending to $c$. Every finite union of the
sets $U_j$ is closed and disjoint from $C$. Hence the indices $j(n)$ for which
$x_n\in U_{j(n)}$ tend to infinity. The estimate just proved yields $c_{j(n)}\to c$,
which gives the last assertion.
\end{proof}

\begin{theorem}
Let $M$ and $M'$ be compact metric spaces, each having infinitely many components and
at least one nondegenerate component. In each space suppose that
\begin{enumerate}
\item[(1)] the nondegenerate components form a null family and have empty interior,
\item[(2)] for each homeomorphism type of nondegenerate component, the union of the
components of that type is dense.
\end{enumerate}
Then $M$ and $M'$ are homeomorphic if and only if $\mathcal{T}(M)=\mathcal{T}(M')$.
\end{theorem}

\begin{proof}
Necessity follows because homeomorphisms preserve components. We prove sufficiency. Fix
compatible metrics on $M$, $M'$ and on $Q_M$, $Q_{M'}$. Enumerate all nondegenerate
components in each space. These families are countable by nullity. By condition (2),
every nonempty clopen set contains a component of each type. We construct paired clopen
sets, starting with $(M,M')$ at level zero.

Suppose that nonempty clopen sets $U\subset M$ and $V\subset M'$ have been paired at
level $n$. At an even level choose the first nondegenerate component $C\subset U$ in the
enumeration for $M$ and a component $D\subset V$ of the same type. At an odd level
choose the first $D\subset V$ in the enumeration for $M'$ and then a component $C\subset
U$ of its type. Fix a homeomorphism $h_C\colon C\to D$. Apply Lemma 6.1 in $U$ and $V$
with $\delta=2^{-(n+1)}$. Denote the resulting sets and points by $(U_j,c_j)$ and
$(V_k,d_k)$, respectively. There exists a permutation $\sigma$ of the positive integers
such that
\[
\lim_{j\to\infty}d_{M'}\bigl(h_C(c_j),d_{\sigma(j)}\bigr)=0 .
\]
Indeed, at odd step $r$ take the first unused $j$ and choose an unused $k$ for which
$d_{M'}(h_C(c_j),d_k)<1/r$. At even step $r$ take the first unused $k$ and choose an
unused $j$ satisfying the same inequality. The last assertion of the lemma, applied on
both sides, ensures that each choice is possible. Every index on either side is
eventually used, and only finitely many pairs are chosen before any given step. Pair
$U_j$ with $V_{\sigma(j)}$ at level $n+1$. Perform the same operation in every pair of
that level.

\smallskip
\noindent\emph{Claim.} Every nondegenerate component is selected at a finite level. If a
component of $M$ were never selected, it would lie in a nested sequence of paired clopen
sets. At each even level a different component preceding it in the fixed enumeration
would have to be selected, which is impossible. The same argument at odd levels applies
to $M'$. Along an infinite nested sequence of paired clopen sets, their intersection on
either side is a singleton. Indeed, compactness gives a nonempty intersection and the
diameter bounds in the component quotient put it in a single component. Every one of the
clopen sets contains that entire component. If it were nondegenerate, it could never
have been selected on this sequence, contrary to the preceding paragraph.

We now define $h\colon M\to M'$ by $h=h_C$ on each selected component $C$. On the
remaining points we use the correspondence between the singleton intersections of paired
nested sequences. This gives a bijection which carries every clopen set occurring in the
construction onto its paired set. To prove continuity, first consider a point given by an
infinite nested sequence. The corresponding clopen sets form a neighborhood basis, by
compactness and the singleton intersection property. Their paired sets have the same
property. Now take $x\in C$, where $C$ was selected in a pair $(U,V)$, and let
$x_m\to x$. Eventually $x_m\in U$. Terms in $C$ have $h(x_m)\to h_C(x)$ along that
subsequence. For the remaining terms write $x_m\in U_{j(m)}$. As every finite union of
the $U_j$ is compact and disjoint from $C$, we have $j(m)\to\infty$. The lemma gives
$d_M(x_m,c_{j(m)})\to0$. Uniform continuity of $h_C$ on $C$ implies
$h_C(c_{j(m)})\to h_C(x)$. Using $h(x_m)\in V_{\sigma(j(m))}$ and
$\sigma(j(m))\to\infty$, we obtain
\begin{align*}
d_{M'}(h(x_m),h_C(x))
&\le \sup_{v\in V_{\sigma(j(m))}}d_{M'}\bigl(v,d_{\sigma(j(m))}\bigr)
 + d_{M'}\bigl(d_{\sigma(j(m))},h_C(c_{j(m)})\bigr)\\
&\quad + d_{M'}\bigl(h_C(c_{j(m)}),h_C(x)\bigr)\longrightarrow0 .
\end{align*}
Thus $h$ is a continuous bijection from a compact space to a Hausdorff space and is a
homeomorphism.
\end{proof}

For a minimal set on a regular curve with infinitely many components, the nondegenerate
components form a null family. Each component has empty interior, since the translates
of any nonempty open subset cover the minimal set and compactness gives a finite
subcover. For each component $C$, the union of its translates is dense. Thus the
hypotheses of Theorem 6.2 hold whenever at least one component is nondegenerate. When
$\mathcal{T}(M)$ is finite, \cite[Theorem 0.2]{ref20} identifies $M$ with the dense
amalgam of one representative of each type. Theorem 6.2 also covers countably many
types.

\subsection{Realization of component types.}

\begin{theorem}
Let $(Y_i)_{i\in I}$ be a nonempty finite or countably infinite family of nondegenerate
regular curves. There are a regular curve $X$, an action of the free group $F_\infty$ on
countably infinitely many generators, and a proper minimal set $M\subset X$ such that
the nondegenerate components form one orbit for each $i\in I$. Every component in that
orbit is homeomorphic to $Y_i$ and has trivial stabilizer. All other components are
singletons, and singleton components occur. If all the $Y_i$ are dendrites, then $X$ is
a dendrite.
\end{theorem}

\begin{proof}
Let $B$ be the countable multigraph with vertices $o$ and $i$ for $i\in I$. For each $i$
join $o$ to $i$ by edges $e_{i,n}$, $n\ge1$. One edge for each $i$ gives a spanning
tree, so the group of deck transformations of the combinatorial universal covering tree
$T\to B$ is isomorphic to $F_\infty$. Denote this group by $G$. It acts freely and
transitively on the lifts of each vertex of $B$ and preserves the edge labels.

For each $i$ choose a sequence $(p_{i,n})_{n\ge1}$ in $Y_i$ such that every nonempty
open subset of $Y_i$ contains infinitely many terms. Replace every vertex $v$ of $T$
lying over $i$ by a copy $P_v$ of $Y_i$. For a vertex $v$ over $o$ let $P_v$ be a single
point. Each lifted edge with label $e_{i,n}$ is represented by an interval joining the
corresponding central point to the copy of $p_{i,n}$ in its other vertex space. Call
these intervals \emph{bridges}.

For a finite nonempty subtree $S\subset T$, let $X_S$ consist of its vertex spaces and
bridges. If $S\subset S'$, retract $X_{S'}$ onto $X_S$ by collapsing each part outside
$S$ to the point at which it attaches to $X_S$. These continuous retractions are
compatible. Define $X$ as their inverse limit. The finite subtrees admit a countable
cofinal exhaustion. Thus $X$ is a compact metric continuum, each $X_S$ embeds
canonically in $X$, and their union is dense. Write $r_S$ for the projection onto $X_S$.
The $G$-action extends to $X$ because it permutes the finite subtrees and commutes with
the corresponding retractions.

Root $T$ at a vertex $c$ over $o$. An end of $T$ is represented by an infinite ray
starting at $c$. The points of $X$ outside the union of the $X_S$ correspond exactly to
these ends. A projection that goes beyond a given vertex or bridge determines the next
edge on the path from $c$. If this path is finite, all projections are those of a point
in its last vertex space or bridge. Otherwise the successive edges form a ray, and the
ray determines all projections. Deleting an edge of $T$ gives two half-trees. The part of
$X$ corresponding to either half-tree is attached to the other part by that bridge.

For fixed $v$, the parts reached from $v$ through distinct incident edges, including
their bridges and attachment points, have diameters tending to zero in any compatible
metric on $X$. To see this, take an increasing exhaustion by finite subtrees containing
$v$. For each fixed subtree in the exhaustion, all but finitely many such parts project
to singletons. The diameters of the fibers of these projections tend uniformly to zero,
since the projections separate points of the compact metric space $X$.

\smallskip
\noindent\emph{Claim 1.} The continuum $X$ is regular. Bridge interiors have their usual
interval neighborhoods. At an end, the part beyond an interior point of a bridge on its
ray is an open neighborhood with one boundary point. These parts form a basis, since any
finite subtree is eventually behind the cut. Now take $x\in P_v$. Given a neighborhood
of $x$, take a smaller one of the form $r_S^{-1}(O)$, where $v\in S$ and $O$ is open in
$X_S$. Choose an open neighborhood $U$ of $x$ in $P_v$ with finite boundary and with
$U\subset O\cap P_v$. Shrink $U$ to avoid the finitely many attachment points in $S$
different from $x$. Include the whole part attached at each point of $U$ through an edge
outside $S$. For the finitely many incident edges of $S$ attached at $x$, include only a
short initial segment, excluding its last point and chosen inside $O$. The resulting
neighborhood is open in $X$. Its boundary consists at most of the boundary points of $U$
in $P_v$ and the finitely many cuts on bridges. Indeed, any point off $P_v$ in a part
attached at a boundary point of $U$ can itself be separated from $P_v$ by a cut on the
initial bridge. This gives the required basis at $x$.

Let $M$ be obtained from $X$ by removing all open bridge interiors. It is closed, proper
and $G$-invariant. Its components are the $P_v$ and the singleton ends. Indeed, any two
such sets are separated by an interior point of a bridge on the tree path between them.
Each $P_v$ is connected. The union of the $P_v$ is dense in $M$, since vertices on a ray
converge to its end.

\smallskip
\noindent\emph{Claim 2.} Every orbit in $M$ meets each side of every bridge. Let $D$ be
one of the half-trees obtained by deleting the corresponding edge of $T$. For an end or
vertex $w$ of $T$ other than $c$, consider the first edge label $e$ on the path from $c$
to $w$. Choose a central vertex $u\in D$ for which the first edge from $u$ toward the
complement of $D$ has label different from $e$. For this, start at any central vertex of
$D$, move away from the deleted edge to an adjacent vertex of some type $i$, and then to
a new central vertex. The last move has infinitely many available labels, so its label
can avoid $e$. The deck transformation carrying $c$ to $u$ preserves labels and hence
sends $w$ into $D$. When $w=c$, transitivity on central vertices suffices. This proves
the assertion also for points of vertex spaces, by applying it to their vertices.

Fix $v$ of type $i$ and $y\in P_v$. Choose distinct $n_j$ with $p_{i,n_j}\to y$, using
the coordinates in $P_v$. At $v$ the bridges labeled $e_{i,n_j}$ enter distinct parts
whose diameters tend to zero. For each $z\in M$, Claim 2 gives $g_j\in G$ carrying $z$
into the part beyond the corresponding bridge. Thus $g_jz\to y$. At a central vertex the
same argument applies to any sequence of distinct incident bridges, all attached at the
central point. It follows that the orbit closure of every $z\in M$ contains every $P_v$
and hence all of $M$. Freeness and transitivity on lifts of a fixed vertex give the
claimed component orbits and stabilizers.

Finally suppose that every $Y_i$ is a dendrite. A simple closed curve in $X$ cannot
intersect an open bridge interior, since a point there separates the two sides of the
bridge and has an interval neighborhood. It would therefore be contained in one component
of $M$, which is impossible. Since $X$ is locally connected, it is a dendrite.
\end{proof}

\subsection{A finitely generated example.}

\begin{proposition}
There is a finitely generated group acting on a dendrite with a proper minimal set whose
nondegenerate components are arcs and whose remaining components are singletons.
\end{proposition}

\begin{proof}
Let $A=\mathbb{Z}+\sqrt2\,\mathbb{Z}$ and $G=A*\langle s\mid s^2=1\rangle$. Consider the
graph whose vertices are $gA$ and whose directed edges, indexed by $g\in G$, go from
$gA$ to $gsA$. The reverse edge is indexed by $gs$. The normal form in the free product
shows this is a tree: a reduced closed path would give a word with alternating
nonidentity $A$-letters and $s$-letters that belongs to $A$. At each vertex place a copy
of the two-end compactification $I=[-\infty,+\infty]$ of the real line. More precisely,
take the vertex spaces in $G\times_A I$, with $(ga,t)$ identified with $(g,t+a)$. The
bridge for the edge indexed by $g$ joins $[g,0]$ to $[gs,0]$. At the vertex $A$, which we
take as the root, the incident bridges attach precisely at the points of $A$, dense in
$I$. Use the inverse limit of the finite subtrees as in Theorem 6.3. The same
neighborhood argument gives a dendrite $X$ and an action of $G$ on it. Removing the open
bridges gives a closed set $M$ whose components are the vertex arcs and the ends.

Fix a point $x\in M$ outside the root arc. It belongs to the part attached at some
$a\in A$. For any $t\in I$ choose distinct $b_j\in A$ so that $a+b_j\to t$. The points
$b_jx$ lie in distinct parts of diameters tending to zero, attached at $a+b_j$. Hence
$b_jx\to t$. If $x$ is on the root arc, first apply $s$ to send it to a different vertex
arc. The orbit closure of every $x\in M$ thus contains the root arc and, by transitivity
on vertices, every vertex arc. Their union is dense in $M$. Finally $A$ is generated by
$1,\sqrt2$, so $G$ is finitely generated.
\end{proof}

\subsection{Triod and circle components.}
In Theorem 6.3 take $Y_1$ to be a triod (three copies of $[0,1]$ with a common endpoint)
and $Y_2$ a circle. For the triod, enumerate the points with rational coordinates in the
interiors of the three arms. For the circle enumerate the rational angles modulo one.
These choices specify all bridge attachments in the construction. The resulting minimal
set has exactly two homeomorphism types of nondegenerate components. Its singleton
components include the central vertices and the ends. Every triod and every circle has
trivial stabilizer, although the union of the copies of either type is dense in the
entire minimal set. For each $m\ge3$, an $m$-od consists of $m$ copies of $[0,1]$ with
one common endpoint. Taking one $m$-od for each $m$ realizes countably many types,
distinguished by the order of the common endpoint.

\subsection{An action of the two-generator free group on a dendrite.}
We use the following special case of the Bass--Serre construction \cite[Ch.~I, \S5.3,
Theorem 12]{ref17}.

\begin{lemma}
Let $H=F(u,v)$ and
\[
\Gamma=\langle u,v,t\mid tut^{-1}=v\rangle .
\]
Then $H$ embeds in $\Gamma$, and $\Gamma$ is freely generated by $u,t$. The graph with
vertex set $\Gamma/H$ and directed edges $g\langle u\rangle$ from $gH$ to $gt^{-1}H$ is a
tree. Left translation defines an action of $\Gamma$ without inversions, transitive on
vertices and on directed edges. The stabilizer $H$ of $H\in\Gamma/H$ has two orbits on
the incident edges, indexed by $H/\langle u\rangle$ and $H/\langle v\rangle$,
respectively. Both orbits are infinite.
\end{lemma}

\begin{proposition}
The free group $F_2$ admits an action on a dendrite with a proper minimal set whose
nondegenerate components are arcs forming a single infinite orbit. The remaining
components are singletons.
\end{proposition}

\begin{proof}
Let $I=[-\infty,+\infty]$ have its order topology. On $\mathbb{R}$ define $S(x)=x+1$ and
\[
A(x)=\begin{cases}
x, & x\le 0,\\
2x, & 0\le x\le 1,\\
x+1, & x\ge 1,
\end{cases}
\]
The map $B=S\circ A^{-1}$ is given by
\[
B(x)=(S\circ A^{-1})(x)=\begin{cases}
x+1, & x\le 0,\\
1+x/2, & 0\le x\le 2,\\
x, & x\ge 2.
\end{cases}
\]
These maps extend to homeomorphisms of $I$ fixing its two endpoints. The map $A$ fixes
$0$, and $B$ fixes $2$. The orbit of $0$ under $\langle A,B\rangle$ is dense in $I$.
Indeed, $BA=S$, and on $[0,1]$ the maps $A^{-1}$ and $A^{-1}S$ are $x\mapsto x/2$ and
$x\mapsto(x+1)/2$. Their iterates starting at $0$ give the dyadic points in $[0,1)$, and
integer translates of these points are dense in $I$. Since $2=S^2(0)$, the two marked
points have the same orbit.

Let $H=F(u,v)$ act on $I$ by $u\mapsto A$ and $v\mapsto B$. Form the HNN extension
\[
\Gamma=\langle u,v,t\mid tut^{-1}=v\rangle .
\]
By Lemma 6.5, $\Gamma=F(u,t)$ and its Bass--Serre tree $T$ has vertices $\Gamma/H$. The
directed edge $g\langle u\rangle$ goes from $gH$ to $gt^{-1}H$. Each vertex has
infinitely many incoming and infinitely many outgoing edges.

Place a copy of $I$ at each vertex. To keep the coordinates equivariant, write its points
as $[g,x]$ with $[gh,x]=[g,hx]$ for $h\in H$. Join $[g,0]$ to $[gt^{-1},2]$ by the bridge
indexed by $g\langle u\rangle$. Replacing $g$ by $gu^m$ fixes both attachment points,
since $u^m$ fixes $0$, $u^mt^{-1}=t^{-1}v^m$, and $v^m$ fixes $2$. Left translation on
the first coordinate acts on this tree of intervals and on the bridges, preserving their
normalized parameters. At $H$ the outgoing bridges attach at the points $[1,h0]$ for
$h\in H$, and the incoming bridges attach at $[1,h2]$. These sets are dense in the
vertex interval. The same holds at every vertex.

Apply the inverse-limit construction over finite subtrees in the proof of Theorem 6.3.
It gives a compact dendrite $X$ and extends the $\Gamma$-action to its ends. Remove the
open bridge interiors to obtain $M$. Its components are the vertex intervals and the
singleton ends. At a fixed vertex the distinct parts reached through its incident bridges
form a null family when each part includes its attachment point.

To prove minimality, first cut any edge of the tree. Every vertex orbit meets both
resulting half-trees by vertex transitivity. Every end orbit also meets both. Indeed,
$\Gamma$ fixes no end, since otherwise $H$ would fix the first edge from $H$ toward that
end, whereas each of its two orbits on incident edges is infinite. For any end $\xi$,
the union of the bi-infinite paths joining distinct members of $\Gamma\xi$ is a nonempty
invariant subtree. By vertex transitivity it is the whole tree. Thus the given edge lies
on a path joining two members of $\Gamma\xi$ on opposite sides. Now fix $z\in M$ and a
point $y$ in a vertex interval. There are distinct incident bridges with attachment
points converging to $y$. The orbit of $z$ enters the far side of each of these bridges.
Since the corresponding parts form a null family, $y\in\overline{\Gamma z}$. The union of
the vertex intervals is dense in $M$, which proves minimality.
\end{proof}

\begin{proof}[Proof of Theorem C]
Proposition 2.7 gives the assertions about the component quotient and its fibers. If all
components are singletons, $M$ is a Cantor set. Otherwise the nullity, empty interiors
and density of each component type, verified after Theorem 6.2, allow us to apply that
theorem. Thus the set $\mathcal{T}(M)$ determines the homeomorphism type of $M$. Theorem
6.3 realizes every prescribed nonempty finite or countable set of types, with trivial
component stabilizers and with a dendrite as ambient space when all prescribed components
are dendrites.
\end{proof}
\section{The free subgroup alternative}

We prove the free-subgroup assertion in Theorem D, and then establish the obstruction
for wreath products. We use the notation and definitions of proximality from Section 2.
We shall use the following results.

\begin{lemma}[{\cite[Lemma 1.5]{ref18}}]
Let a countable group $G$ act minimally on an infinite compact metric space $Y$. For
every countable set $C\subset Y$ and every nonempty finite set $F\subset Y$, there are
$g_n\in G$ and a finite set $K\subset Y$ such that $g_nF\to K$ in the Hausdorff metric
and $K\cap C=\emptyset$.
\end{lemma}

\begin{theorem}[{\cite[Theorem 3.25]{ref16}}]
Let $K$ be a compact convex subset of a locally convex space, and let $S\subset K$. If
$K$ is the closed convex hull of $S$, then $\ext K\subset S$.
\end{theorem}

\begin{theorem}[{\cite[Chapter III, Theorem 2.3]{ref8}}]
Let $K$ be a nonempty compact metrizable convex subset of a locally convex space, with a
continuous affine $G$-action. Suppose that $K$ contains no proper nonempty closed convex
invariant subset. Then $\ext K$ is minimal and strongly proximal.
\end{theorem}

\begin{proof}
Put $Y=\ext K$. For every $\theta\in K$, the closed convex hull of $G\theta$ is $K$.
Theorem 7.2 gives $\ext K\subset\overline{G\theta}$. Taking $\theta\in Y$ proves
minimality. Let $\lambda\in\mathcal{P}(Y)$, and let $b(\lambda)\in K$ be its barycenter.
Fix $\theta_0\in\ext K$ and choose $g_n\in G$ such that $g_nb(\lambda)\to\theta_0$. By
continuity and equivariance of the barycenter map, every limit of $g_{n*}\lambda$ has
barycenter $\theta_0$. Such a probability is $\delta_{\theta_0}$. Indeed, otherwise a
continuous affine function nonconstant on its support splits it into two probabilities
with distinct barycenters, contrary to extremality. Thus
$g_{n*}\lambda\to\delta_{\theta_0}$.
\end{proof}

\begin{lemma}
Let $G$ act on a regular curve $X$, and let $L\subset X$ be a nonempty closed invariant
set. If $Y\subset\mathcal{P}(L)$ is compact, invariant, and consists of atomless
measures, then the action on $Y$ is equicontinuous.
\end{lemma}

\begin{proof}
Given $g_n\in G$, use Lemma 2.2 to pass to a subsequence such that the graphs of $g_n|_L$
converge in the Hausdorff metric to a compact relation $R\subset L^2$. Its first
projection is $L$. Write $R_x=\{y:(x,y)\in R\}$.

\smallskip
\noindent\emph{Claim.} Only countably many fibers of $R$ are nonsingletons. Fix
$\varepsilon>0$. Choose distinct $x_1,\dots,x_k$ with $\diam R_{x_i}\ge\varepsilon$ and
pairwise disjoint connected open neighborhoods $U_i$ of $x_i$ in $X$. Two points of each
$R_{x_i}$ are at distance at least $\varepsilon$. For all large $n$, they are
approximated by graph points with first coordinates in $U_i$ and second coordinates at
distance at least $\varepsilon/2$. The disjoint connected sets $g_nU_i$ therefore have
diameter at least $\varepsilon/2$. Lemma 2.3 bounds $k$ independently of the chosen
points. There are consequently only finitely many fibers of diameter at least
$\varepsilon$, proving the claim.

Let $D$ be the countable set in the claim. Suppose that $\nu_n,\nu'_n\in Y$ tend to the
same $\nu\in Y$. Pass to further subsequences such that $(\mathrm{id},g_n)_*\nu_n$ and
$(\mathrm{id},g_n)_*\nu'_n$ have weak limits on $L^2$. Both limits are supported on $R$
and have first marginal $\nu$. Since $\nu(D)=0$, both are concentrated on the graph over
$L\setminus D$. The function defined by this graph is continuous on that subspace, by
closedness of $R$ and compactness of $L$. Hence the two limits coincide, and so do their
second marginals. If the action on $Y$ were not equicontinuous, compactness would provide
$\nu_n,\nu'_n$ tending to a common limit and $g_n\in G$ whose images stay a fixed
positive distance apart. The preceding argument gives a contradiction.
\end{proof}

\begin{theorem}
Let $G$ act on a regular curve $X$, let $L\subset X$ be a minimal set, and let
$\pi\colon L\to Y$ be a factor map onto an infinite compact metric system. Suppose that
some nonempty finite set $A\subset Y$ satisfies the following condition. For every finite
$E\subset Y$ and every open neighborhood $W$ of $A$, there is $g\in G$ such that
$gE\subset W$. Then the image of $G$ in $\Homeo(Y)$ contains a nonabelian free subgroup.
If $Y$ is proximal, it is extremely proximal.
\end{theorem}

\begin{proof}
\emph{Claim.} If open sets $O,V\subset Y$ each contain a translate of $A$, then
$s(Y\setminus V)\subset O$ for some $s\in G$. Write $tA\subset O$. Choose an open set
$U\subset X$ with finite boundary such that
\[
\pi^{-1}(tA)\subset U,\qquad U\cap L\subset\pi^{-1}(O).
\]
A finite union of basic neighborhoods covering $\pi^{-1}(tA)$ has these properties. Let
$D\subset L$ be countable and dense. The hypothesis gives $g_n\in G$ with
$\dist(g_n\pi(x),tA)\to0$ for every $x\in D$. Thus $g_nx\in U$ eventually for each
$x\in D$. Passing to a subsequence using Lemma 2.2, assume that $g_n^{-1}\partial_X U$
converges in the Hausdorff metric to a finite set $F$. If $\partial_X U=\emptyset$, take
$F=\emptyset$. For every compact $B\subset L\setminus F$, we have $g_nB\subset U$ for all
large $n$. To see this, cover $B$ by finitely many connected open subsets of $X$ whose
closures miss $F$ and which meet $L$. Each meets $D$. Eventually its image meets $U$ and
misses $\partial_X U$, so lies in $U$. Choose $h\in G$ such that $h\pi(F\cap L)\subset V$.
The compact set
\[
B=L\cap\pi^{-1}h^{-1}(Y\setminus V)
\]
misses $F$. Hence $s=g_nh^{-1}$ proves the claim for all large $n$.

Lemma 7.1, applied successively with $C$ equal to the union of the preceding translates,
gives four pairwise disjoint translates of $A$. Choose pairwise disjoint open
neighborhoods $U_+,U_-,V_+,V_-$ of these translates whose union misses a point $y\in Y$.
By the claim there are $s,t\in G$ with $s(Y\setminus U_-)\subset U_+$ and
$t(Y\setminus V_-)\subset V_+$. Every nonempty reduced word in $s,t$ sends $y$ into the
set corresponding to its first letter. Their images therefore freely generate a subgroup
of $\Homeo(Y)$, and $s,t$ freely generate a subgroup of $G$.

Suppose now that $Y$ is proximal. Every finite tuple can be moved into any prescribed
nonempty open set. Indeed, a minimal subset of its orbit closure projects on each pair of
coordinates to a minimal subset of $Y^2$. Proximality makes each such projection meet,
and therefore lie in, the diagonal. The minimal subset of the tuple orbit closure is thus
the full diagonal. We may consequently take $A$ to be a singleton. Its translates meet
every nonempty open set, so the claim applies to arbitrary nonempty open sets $O,V$. This
proves extreme proximality.
\end{proof}

\begin{theorem}
Let a countable group $G$ act on a regular curve $X$. If a minimal set $M\subset X$
admits no $G$-invariant probability measure, then the image of $G$ in $\Homeo(M)$
contains a nonabelian free subgroup. In particular, if $X$ admits no invariant
probability measure, then $G$ contains a nonabelian free subgroup.
\end{theorem}

\begin{proof}
By compactness, choose a minimal nonempty compact convex invariant subset
$K\subset\mathcal{P}(M)$. Put $Y=\ext K$. By Theorem 7.3, $Y$ is minimal and strongly
proximal. It has more than one point, since a singleton would give an invariant
probability on $M$.

\smallskip
\noindent\emph{Claim 1.} There is $a>0$ such that every $\theta\in Y$ has maximal atom
mass $a$. Some $\theta\in Y$ has an atom. Otherwise Lemma 7.4 would make $Y$
equicontinuous. Two distinct points of an equicontinuous system cannot be proximal,
whereas strong proximality contracts their equal-weight average to a point mass. The
function
\[
a(\theta)=\sup_{x\in M}\theta(\{x\})
\]
is invariant and upper semicontinuous. Indeed, if $\theta_n\to\theta$ have atoms of mass
at least $c>0$ at $x_n\to x$, every closed neighborhood of $x$ has $\theta$-mass at least
$c$, and hence $\theta(\{x\})\ge c$. Minimality now makes $a(\theta)$ a positive
constant.

\smallskip
\noindent\emph{Claim 2.} The sets $A(\theta)=\{x\in M:\theta(\{x\})=a\}$ have constant
finite cardinality and depend continuously on $\theta$ in the Hausdorff topology. Each
$A(\theta)$ is nonempty and has at most $1/a$ points. For every integer $j$, the set
$\{\theta\in Y:|A(\theta)|\ge j\}$ is closed. Indeed, along $\theta_n\to\theta$, take
limits of $j$ chosen atoms. Each limit has mass $a$, and no two limits coincide, since
their common limit would have mass at least $2a$. Invariance and minimality therefore
give $|A(\theta)|=k$ for a fixed $k$. The same argument shows that every subsequential
limit of $A(\theta_n)$ consists of exactly the $k$ maximal atoms of $\theta$. This proves
the claim.

The map $\theta\mapsto A(\theta)$ is equivariant. Its union is compact and invariant, so
equals $M$. Fix $\theta_0\in Y$ and put $A=A(\theta_0)$. For a nonempty finite set
$E\subset M$, choose $\theta_x\in Y$ with $x\in A(\theta_x)$ for each $x\in E$. Strong
proximality and minimality give a sequence in $G$ sending the uniform average of the
$\delta_{\theta_x}$ to $\delta_{\theta_0}$. The same sequence sends every $\theta_x$ to
$\theta_0$. By Claim 2, every neighborhood $W$ of $A$ therefore contains $gE$ for some
$g\in G$. The set $M$ is infinite, since a finite orbit carries its uniform invariant
probability. Theorem 7.5, applied to the identity factor of $M$, gives the required free
subgroup. If $X$ has no invariant probability, neither does any minimal set. Lifts of
free generators of the image freely generate a subgroup of $G$.
\end{proof}

\subsection{Commuting conjugates and wreath products.}

\begin{lemma}
Let $H$ act minimally on a closed subset $N$ of a regular curve by restrictions of
homeomorphisms of the curve. If $N$ admits no $H$-invariant probability measure, then
the centralizer of this action in $\Homeo(N)$ is finite.
\end{lemma}

\begin{proof}
The proof of Theorem 7.6 gives a nonempty finite set $A\subset N$ such that every finite
subset of $N$ can be moved into every neighborhood of $A$. Let $P$ be the closure of the
acting maps in the compact product $N^N$. The closed sets
\[
\{p\in P:p(x)\in A\},\qquad x\in N,
\]
have the finite intersection property. For finitely many prescribed points, move them
into a decreasing neighborhood basis of $A$ and take a cluster point in $P$. Thus some
$p\in P$ satisfies $p(N)\subset A$. Every homeomorphism $c$ commuting with $H$ commutes
with $p$, by pointwise approximation and continuity. It therefore permutes the nonempty
finite set $p(N)$. A homeomorphism commuting with $H$ and fixing one point fixes its
dense $H$-orbit. Evaluation at any point of $p(N)$ consequently injects the centralizer
into $p(N)$.
\end{proof}

\begin{proposition}
Let $G$ be countable, and let $B\triangleleft G$ with $G/B$ amenable. Suppose that every
finitely generated subgroup $H\le B$ has a conjugate $sHs^{-1}$ commuting elementwise with
$H$. Then every nonempty closed invariant subset of every regular-curve $G$-action
supports a $G$-invariant probability measure.
\end{proposition}

\begin{proof}
Suppose that such a subset $L$ has no $G$-invariant probability. It has no $B$-invariant
probability either. Otherwise $G/B$ would act affinely on the nonempty compact convex set
of $B$-invariant probabilities on $L$. Lemma 2.1 would give a $G$-invariant member.
Compactness of $\mathcal{P}(L)$ gives a finite subset of $B$ with no common invariant
probability. The subgroup $H$ it generates has no invariant probability on $L$. Choose
$H'=sHs^{-1}$ commuting with $H$. Since $sL=L$, neither does $H'$ have an invariant
probability on $L$.

Let $N\subset L$ be $H$-minimal. The group $H$ has no finite orbit in the ambient curve.
Otherwise Corollary 3.3 would make $N$ finite or equicontinuous Cantor, and in either
case it would carry an invariant probability. By the same corollary, $N$ is the unique
ambient $H$-minimal set. Since $H'$ commutes with $H$, it preserves $N$ and acts through
the centralizer of $H|_N$. This centralizer is finite by Lemma 7.7. The uniform
probability on any $H'$-orbit in $N$ is then an $H'$-invariant probability on $L$, a
contradiction.
\end{proof}

\begin{corollary}
Let $D$ be countable, and let $T$ be infinite, countable and amenable. For the restricted
wreath product
\[
D\wr T=\Bigl(\bigoplus_{t\in T}D_t\Bigr)\rtimes T,\qquad D_t\cong D,
\]
where $T$ acts on the indices by left translation, every minimal set of every
regular-curve action supports an invariant probability measure. In particular,
$F_2\wr\mathbb{Z}$ is finitely generated and contains $F_2$, but admits no minimal set
with infinitely many components and a nondegenerate component.
\end{corollary}

\begin{proof}
A finitely generated subgroup of the base is supported in a finite set $J\subset T$.
Choose $s\in T\setminus JJ^{-1}$. Then $J\cap sJ=\emptyset$, so the corresponding
coordinate subgroups commute. Apply Proposition 7.8. For $D=F_2$ and $T=\mathbb{Z}$, a
copy of $F_2$ at the origin and the unit shift generate the wreath product. Theorem 4.4
excludes the stated minimal sets.
\end{proof}

\begin{proof}[Proof of Theorem D]
A minimal set with infinitely many components and a nondegenerate component has no
invariant probability by Theorem 4.4. The first assertion follows from Theorem 7.6. The
second is Proposition 6.6, and the third is Corollary 7.9.
\end{proof}
\section{Consequences for a prescribed group}

By Theorem A, a minimal set with finitely many components is connected or a finite
orbit. For infinitely many components, Theorem 6.2 determines the space from its
nondegenerate component types, and Subsection 6.2 realizes every finite or countable set
of such types when the group may vary. We record two restrictions for a prescribed
group.

\subsection{Expansion on the component quotient.}
A subgroup $H\le G$ is coamenable if the left action of $G$ on $G/H$ admits an invariant
mean.

\begin{proposition}
Let $M$ be a minimal set of a regular-curve $G$-action with infinitely many components,
at least one nondegenerate, and let $q\colon M\to Q$ be the component quotient. There are
a finite symmetric set $S\subset G$ and $\eta>0$ such that
\[
|SA\setminus A|\ge\eta|A|
\]
for every finite nonempty $A\subset Q$. In particular, no component stabilizer is
coamenable in $G$.
\end{proposition}

\begin{proof}
By Theorem 4.4, $M$ has no invariant probability, and by Proposition 2.7, neither does
$Q$. If the bound failed, choose increasing finite symmetric sets $S_n$ exhausting $G$
and finite nonempty $A_n\subset Q$ with $|S_nA_n\setminus A_n|<|A_n|/n$. For $s\in S_n$,
$|sA_n\triangle A_n|<2|A_n|/n$. A weak subsequential limit of the uniform probabilities
on $A_n$ would be invariant, a contradiction. If a component stabilizer $H$ were
coamenable, an invariant mean on $G/H$ applied to the restrictions of continuous
functions on $Q$ to that orbit would give an invariant probability on $Q$.
\end{proof}

\subsection{Finite orbits and realization.}

\begin{corollary}
Suppose that a countable group $G$ contains no nonabelian free subgroup. Every
regular-curve $G$-action has a finite orbit if and only if there is no homomorphism
$G\to\Isom(S^1)$ with infinite image.
\end{corollary}

\begin{proof}
An infinite subgroup of $\Isom(S^1)$ has a dense rotation subgroup and acts without
finite orbits. Conversely, suppose an action has no finite orbit. Its unique minimal set
supports an invariant probability by Theorem 7.6. An equicontinuous Cantor minimal set
has a nontrivial finite factor, which would give an ambient finite orbit by Theorem 3.4.
The minimal set is therefore a circle isometry action or a split extension of one. In
either case the circle action gives a homomorphism with infinite image.
\end{proof}

For minimal sets with infinitely many components and at least one nondegenerate
component, we have
\[
G\twoheadrightarrow F_2
\ \Longrightarrow\ G \text{ admits such a minimal set}
\ \Longrightarrow\ F_2\le G .
\]
The first implication follows by pulling back the action in Proposition 6.6. Its converse
fails for $(\mathbb{Z}^2)*C_2$ by Proposition 6.4. Every homomorphism from this group to
$F_2$ kills the involution, so has abelian image. The second implication is Theorem D,
and its converse fails for $F_2\wr\mathbb{Z}$ by Corollary 7.9.

\begin{question}
Which intrinsic group-theoretic condition characterizes countable groups admitting a
minimal set with infinitely many components, at least one nondegenerate, on a regular
curve?
\end{question}

The induced action on components must also satisfy Proposition 8.1. The preceding
examples show that neither containing a nonabelian free subgroup nor admitting a
nonabelian free quotient gives a characterization.

\end{document}